\documentclass[pre,twocolumn,showkeys,showpacs,floatfix]{revtex4-1}
\usepackage{epsfig,graphicx,subfigure,url}
\graphicspath{{Figures/}}
\usepackage{amsmath,amsfonts}
\usepackage{mathptmx,paralist}
\pdfoutput=1
\usepackage{url,color}
\usepackage[pdftex,colorlinks]{hyperref}
\definecolor{darkblue}{cmyk}{1,0,0,0.8}
\definecolor{darkred}{cmyk}{0,1,0,0.7}
\hypersetup{anchorcolor=black, citecolor=darkblue,
  filecolor=darkblue,
  menucolor=darkblue,pagecolor=darkblue,urlcolor=darkblue,linkcolor=darkblue} 
\usepackage{microtype}
\renewcommand{\epsilon}{\varepsilon}
\newcommand{\R}{\mathbb{R}}

\renewcommand{\d}{\mathrm{d}}

\newcommand{\fref}[1]{Fig.~\ref{#1}}
\newcommand{\Fref}[1]{Figure~\ref{#1}}

\makeatletter
\def\diff{\@ifnextchar[{\@diffwith}{\@diffwithout}}
\def\@diffwith[#1]#2#3{\frac{\md^{#1} #2}{\md #3^{#1}}}
\def\@diffwithout#1#2{\frac{\md #1}{\md #2}}
\makeatother
\begin{document}
\title{Systematic experimental exploration of bifurcations with
  non-invasive control} 
\author{D.A.W.\ Barton} 
\affiliation{Department of Engineering Mathematics, University of
  Bristol, Queen's Building, University Walk, Bristol, BS8 1TR, U.K.}
\author{J.\ Sieber}
\affiliation{College of Engineering, Mathematics and Physical
  Sciences, University of Exeter, Exeter, EX4 4QF, U.K.}
\begin{abstract}
  We present a general method for systematically investigating the
  dynamics and bifurcations of a physical nonlinear experiment. In
  particular, we show how the odd-number limitation inherent in
  popular non-invasive control schemes, such as (Pyragas) time-delayed
  or washout-filtered feedback control, can be overcome for tracking
  equilibria or forced periodic orbits in experiments. To demonstrate
  the use of our non-invasive control, we trace out experimentally the
  resonance surface of a periodically forced mechanical nonlinear
  oscillator near the onset of instability, around two saddle-node
  bifurcations (folds) and a cusp bifurcation.
\end{abstract}
\keywords{experimental bifurcation analysis, non-invasive control,
  odd-number limitation} \pacs{05.45.Gg,45.80.+r,02.30.Oz}
\maketitle

\noindent Feedback control is not only of interest as a tool for
system manipulation in the classical control engineering sense, but it
can also be used for model verification and discovery, if one can
ensure that it is \emph{non-invasive}. Washout-filtered feedback and
Pyragas' time-delayed feedback (TDF) are popular feedback control
schemes that are automatically non-invasive \cite{AWC94,P92}. For
example, TDF feeds a signal $u(t)=k^T[x(t)-x(t-T)]$ back into the
experimental dynamical system, where $x$ is some system output
(possibly processed), $k$ is a vector or matrix of control gains and
$T$ is an a-priori chosen time delay. If the time delay equals the
period of a periodic orbit $x(t)$ ($t\in[0,T]$) of the uncontrolled
dynamical system and the gains are such that the control is
stabilizing, then the controlled system will also have the periodic
orbit $x$, because the control input $u$ vanishes for all time (that
is, the control becomes non-invasive). However, the control has
changed the stability of $x$, making it visible in the experiment.

While sometimes TDF (or its extended version \cite{GSCS94}) is used
for engineering purposes (suppression of period doublings leading to
chaos \cite{YKYMH09,AP04}), non-invasiveness is not essential in these
applications. Thus, in these cases the delayed term $x(t-T)$ in the
feedback loop could have been replaced by a periodic reference signal
$x^*(t)$ approximately resembling the desired behaviour. TDF
draws interest mostly in the scientific community because it enables
experimenters to explore dynamical phenomena such as equilibria and
periodic orbits of the original uncontrolled system regardless of
their dynamical stability
\cite{KBPOMBRE01,SWH11,SHWSH06,LWP01,CCL96,LBJ10}.

Systematic studies that try to explore the parameter space of the
uncontrolled system and that use TDF or washout filters to stabilize
equilibria and periodic orbits non-invasively encounter a major
difficulty: it is not known under which conditions one can find
control gains $k$ that successfully stabilize a periodic orbit
\cite{SS07}. This is in contrast to classical feedback control where
one feeds $u(t)=k^T[x(t)-x^*(t)]$ with a pre-determined reference
signal $x^*(t)$ back into the system. For classical feedback control
it is known that, if $x^*$ corresponds to an equilibrium or periodic
orbit of the uncontrolled system then, under some genericity
assumptions (controllability and observability), one can always
locally stabilize the equilibrium or periodic orbit even if it has
arbitrarily many unstable eigenvalues \cite{S98}. Experimental and
theoretical studies have explored the limits of applicability of TDF
and restrictions on the gains due to instabilities caused by the TDF
feedback term $k^T[x(t)-x(t-T)]$ \cite{LBJ10,FFHS10,FS11,FFGHS07}. In
particular, there are topological restrictions (the \emph{odd-number
  limitation} \cite{N04,HA12}) which guarantee that TDF cannot
possibly work in some of the most common cases. One common scenario
where it would be natural to use TDF is ruled out by the odd-number
limitation: an equilibrium of an autonomous system or a periodic orbit
of a periodically forced system in the vicinity of a system parameter
setting where it makes a fold (called saddle-node bifurcation, see
\fref{fig:simple-evolution} for a typical bifurcation diagram). Even
the unstable controller proposed in \cite{TMPP07} fails to stabilize
uniformly near the fold.

In this paper, we present a simple alternative approach to exploring
bifurcation scenarios, including the unstable branches. Our approach
exploits the fact that the goal of the experiment is a parameter study
in a system parameter $p$, rather than finding a single equilibrium or
periodic orbit at a specified parameter value. It is applicable
whenever the feedback control is achieved by (effectively) varying the
same system parameter $p$ that one wants to use for the bifurcation
diagram.

Finally, we illustrate the power of this approach by tracking branches of
periodic orbits in a physical experiment to produce a solution surface
that shows two fold curves meeting at a cusp bifurcation.

\section{Steady-state branch tracking}
\label{sec:method}
\noindent 
We first explain the basic approach applied to tracking a branch of
equilibria in the case of a single-input single-output
system. A stable and an unstable branch, connected by a fold, are
traced out as a function of the system parameter. We then generalize
the approach to the case of tracking periodic orbits and show its
application in an experimental setting. As our method does not rely on
knowledge of the state of the model we do not distinguish between
state and output, calling the output $x$.

\subsection{Equilibria}
\label{sec:method:eq}
\begin{figure}
  \centering
  \includegraphics[width=\columnwidth]{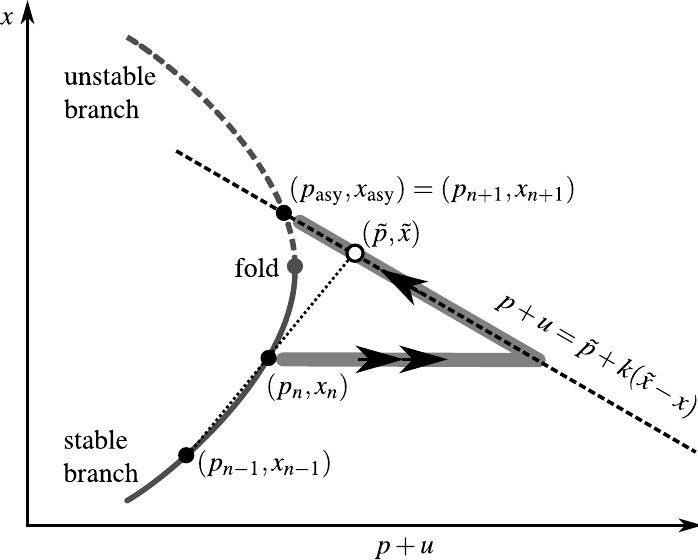}
  \caption{Illustration for tracing out unstable branches near folds
    with feedback control. A family of equilibria (stable=solid dark
    gray, unstable=dashed dark gray) is tracked through a saddle-node
    bifurcation (fold). The thick light gray lines with arrows
    illustrate the dynamical behavior: fast=double arrows,
    slow=single arrow).}
  \label{fig:simple-evolution}
\end{figure}

\noindent Suppose that we have an experiment with a scalar output $x(t)$ and a
scalar system parameter $p$, which has a bifurcation diagram for its
equilibria as illustrated in \fref{fig:simple-evolution}: a stable and
an unstable branch of equilibria meeting in a fold (saddle-node
bifurcation). We assume that we can use the parameter $p$ as a
(scalar) control input such that control $u$ is added to the parameter
$p$.  Consider what happens if we pick a point $(\tilde p,\tilde x)$ in the
$(p,x)$-plane and apply a simple proportional feedback controller of
the form
\begin{equation}
  \label{eq:controlinp}
  p+u(t) = \tilde p + k(\tilde x - x(t)),
\end{equation}
where $k>0$ is the control gain (see \fref{fig:simple-evolution}). The
linear relation \eqref{eq:controlinp} restricts the dynamics of the
experiment to move along the tilted dashed line in
\fref{fig:simple-evolution}. The equilibria of the
feedback-controlled system are the intersections of this tilted line
with the curve of equilibria in the $(p,x)$-plane. In
\fref{fig:simple-evolution} this intersection point is
$(p_\mathrm{asy},x_\mathrm{asy})$. Close to the fold and for
sufficiently small $k$ the underlying dynamics is one-dimensional such
that the intersection point $(p_\mathrm{asy},x_\mathrm{asy})$ corresponds to a stable
equilibrium of the experiment with feedback control
\eqref{eq:controlinp} (and, at the same time, to a possibly unstable
equilibrium of the uncontrolled system).

This simple trick illustrated in \fref{fig:simple-evolution} permits
one to trace out branches of equilibria around folds with a
continuation procedure. Assume that we have already found two
equilibria $(p_{n-1},x_{n-1})$ and $(p_n,x_n)$ along the branch, and
that the experiment is currently at equilibrium $(p_n,x_n)$. First we
use a secant approximation to generate a prediction for the next
equilibrium, which gives
\begin{equation}\label{eq:pred}
  (\tilde p, \tilde x)=(p_n,x_n)+h[(p_n,x_n)-(p_{n-1},x_{n-1})]\text{,}
\end{equation}
where $h=1$ in \fref{fig:simple-evolution}. (The prediction step $h$
can be chosen adaptively to ensure the desired resolution of the
equilibrium branch.)  Then the experiment is run with the feedback
control \eqref{eq:controlinp} based on the point $(\tilde p,\tilde x)$
determined by \eqref{eq:pred}. Once the transients have settled
to a constant value, the next equilibrium along the branch is given by
\begin{align*}
  p_{n+1}&:=p_\mathrm{asy}=\lim_{t\to\infty}\ \tilde p+k(\tilde x-x(t))\mbox{,}\\
  x_{n+1}&:=x_\mathrm{asy}=\lim_{t\to\infty}\ x(t)\mbox{.}
\end{align*}
Then we can repeat the procedure by picking the next prediction using
\eqref{eq:pred} for index $n+1$, finding the next equilibrium along the
branch. 

\Fref{fig:simple-evolution} also illustrates the dynamics of the
system after one sets the feedback control to \eqref{eq:controlinp}
with \eqref{eq:pred}. The system parameter and the control input
adjust immediately such that initially $(p+u(t),x(t))$ jumps rapidly
in horizontal direction (thick light gray line with double arrows in
\fref{fig:simple-evolution}). Then the system follows the dynamics
imposed by the feedback control \eqref{eq:controlinp} along the tilted
line given by \eqref{eq:controlinp}, gradually settling to its
equilibrium (thick light gray line with single arrow in
\fref{fig:simple-evolution}).

The approach illustrated in \fref{fig:simple-evolution} is
non-invasive: every equilibrium of the feedback controlled system
corresponds to an equilibrium of the uncontrolled system and a-priori
knowledge of the equilibrium of the uncontrolled system is not
necessary. 

\subsection{Periodic orbits}
\label{sec:method:po}
\noindent For our demonstration experiment, which is periodically
forced, we need to generalize the approach to enable it to trace
periodic steady states of periodically forced systems with
proportional-plus-derivative (PD) control. Here the input parameter is
harmonic forcing, that is,
\begin{displaymath}
  p(t) = a\cos(\omega t)+b\sin(\omega t)\mbox{.}
 \end{displaymath}
This allows for an arbitrary phase shift in the forcing and simplifies
the method below.  We perform a parameter study in the forcing
amplitude $r=\sqrt{a^2+b^2}$. Mimicking the approach of
\fref{fig:simple-evolution}, we pick a harmonic forcing amplitude
$(a^*,b^*)$
\begin{equation}\label{eq:pstarfcoeffs}
p^*(t)=a^*\cos(\omega t)+b^*\sin(\omega t)
\end{equation}
(playing the role of $\tilde p$ in \fref{fig:simple-evolution}) and an
arbitrary periodic reference signal (expanded to finitely many
Fourier modes)
\begin{equation}\label{eq:ystarfcoeffs}
  x^*(t) = \frac{A^*_0}{2} + \sum_{j=1}^{m}A^*_j\cos(j\omega t) +
  B^*_j\sin(j\omega t)\mbox{,}
\end{equation}
(playing the role of $\tilde x$ in \fref{fig:simple-evolution}), and
apply the PD feedback law (an idealized version of what we do in the
experiment in section~\ref{sec:setup})
\begin{equation}\label{eq:controlinppd}
  p(t)+u(t)=p^*(t)+k_p(x^*(t)-x(t))+k_d(\dot x^*(t)-\dot x(t))\mbox{,}
\end{equation}
where $x(t)$ is the output of the system. Assuming that the PD control
is stabilizing, the system will settle into a periodic steady-state
output (also expanded to Fourier modes):
\begin{equation}\label{eq:yfcoeffs}
  x_\mathrm{asy}(t) = \frac{A_0}{2} + \sum_{j=1}^{m}A_j\cos(j\omega t) +
  B_j\sin(j\omega t)\mbox{.}
\end{equation}
We notice that the experiment with feedback control
\eqref{eq:controlinppd} also has periodic input after the transients
have settled (the right-hand side of \eqref{eq:controlinppd} is
periodic for periodic $x^*$ and $x=x_\mathrm{asy}$). The amplitude of
the forcing at the fundamental frequency $\omega$ equals
$r=\sqrt{a^2+b^2}$, where
\begin{equation}
  \begin{aligned}
    a&=a^* + k_p(A_1^* - A_1) + \omega k_d(B_1^* - B_1)\mbox{,}\\
    b&=b^* + k_p(B_1^* - B_1) + \omega k_d(A_1 - A_1^*)\mbox{.}
  \end{aligned}\label{eq:pupdate}
\end{equation}
However, the forcing is not purely harmonic asymptotically due to the
presence of nonlinearities. Even if
the reference signal $x^*$ is harmonic, the output $x_\mathrm{asy}$
will not be harmonic because of the nonlinearities of the experimental
system. The non-harmonic Fourier coefficients of $u(t)$ are
\begin{align*}
  A^u_0&=k_p(A_0^*-A_0)\mbox{,}\\
  A^u_j&=k_p(A_j^* - A_j) + j\omega k_d(B_j^* - B_j)\mbox{\quad ($j>1$),}\\
  B^u_j&=k_p(B_j^* - B_j) + j\omega k_d(A_j - A_j^*)\mbox{\quad ($j>1$).}
\end{align*}
If these are zero then the forcing $p(t)+u(t)$ will be harmonic with
amplitude $r=\sqrt{a^2+b^2}$ such that the point $(r,x^*)$ will be on
the branch of periodic orbits. The requirement for
$[A^u_0,(A^u_j,B^u_j)_{j=2}^m]$ to be zero is a nonlinear system of
$2m-1$ equations for the non-harmonic Fourier coefficients
$[A^*_0,(A^*_j,B^*_j)_{j=2}^m]$ of the reference signal $x^*$ (which
are $2m-1$ variables). This nonlinear system can be written as a
nonlinear fixed-point problem:
\begin{align}
  \nonumber
  0&=A^u_0=A^u_j=B^u_j \mbox{\quad ($j>1$)\quad if and only if}\\
  \label{eq:zeroproblem}
  X^*&=X(X^*)\mbox{,\quad where}\\[0.5ex]
  \nonumber
  X^*&=[A^*_0,(A^*_j,B^*_j)_{j=2}^m]\mbox{,\ and\ }
  X= [A_0,(A_j,B_j)_{j=2}^m]\mbox{.}
\end{align}
The output $x_\mathrm{asy}(t)$ (and, thus, its vector of non-harmonic
Fourier coefficients $X$) also depends on the harmonic amplitudes
$(a^*,b^*)$ and $(A^*_1,B^*_1)$, which act as parameters in
\eqref{eq:zeroproblem} but were omitted as arguments in
\eqref{eq:zeroproblem}. In general the fixed-point problem
\eqref{eq:zeroproblem} has to be solved with a Newton iteration (this
is what \cite{SGNWK08,Barton2012,BB11,Bureau2011,Bureau2012} do for
all Fourier coefficients). However, in many practical experiments it
may be sufficient to apply a simple fixed point iteration to the
fixed-point problem (this is the approach we take in
section~\ref{sec:setup}):
\begin{equation}
  \label{eq:picardj}
  X^*_{k+1} := X(X^*_k)\mbox{.}
\end{equation}


\paragraph*{Algorithm}
In summary, the procedure to find a new periodic orbit on the branch
looks as follows.  Assume that we have found already two previous
points along the branch of periodic orbits, namely $(p_{n-1},x_{n-1})$
and $(p_n,x_n)$. The inputs $p_{n-1}$ and $p_n$ are harmonic and the
outputs $x_{n-1}$ and $x_n$ are periodic. We use a secant
approximation to generate a prediction for the next solution point,
which gives
\begin{equation}\label{eq:predperiod}
  (\tilde p, \tilde x)=(p_n,x_n)+h[(p_n,x_n)-(p_{n-1},x_{n-1})]
\end{equation}
(again $h=1$ by default, but can be chosen adaptive).  We set
$p^*=\tilde p$ and $x^*=\tilde x$, and repeat the following procedure
until convergence.
\begin{compactenum}
\item Run the experiment with PD feedback law \eqref{eq:controlinppd}
  and $p^*$ and $x^*$ as given. Wait until the experiment settles to a
  periodic output $x_\mathrm{asy}$.
\item Extract the Fourier coefficients:
  \begin{compactitem}
  \item $[A_0^*,(A^*_j,B^*_j)_{j=1}^m]$ from $x^*$ according to
    \eqref{eq:ystarfcoeffs},
  \item $[A_0,(A_j,B_j)_{j=1}^m]$ from $x_\mathrm{asy}$ according to
    \eqref{eq:yfcoeffs}.
  \end{compactitem}
\item Check if the root-mean-square error
  \begin{displaymath}
    \sqrt{(A^*_0-A_0)^2)+\sum_{j=2}^m(A^*_j-A_j)^2+(B^*_j-B_j)^2}
  \end{displaymath}
  is smaller than the desired tolerance. (Note that the index $1$ is
  skipped in the sum!) If yes, finish the iteration.  Otherwise, set
  \begin{displaymath}
    A^*_{0,\mathrm{new}}:=A_0\mbox{,\ } 
    A^*_{j,\mathrm{new}}:=A_j\mbox{,\ }
    B^*_{j,\mathrm{new}}:=B_j\mbox{\ ($j>1$),}
  \end{displaymath}
  set the reference signal $x^*$ to these new non-harmonic Fourier
  coefficients according to \eqref{eq:ystarfcoeffs} (keeping $A^*_1$
  and $B^*_1$ as before), and repeat from step 1.
\end{compactenum}
After the iteration is finished, the accepted point on the branch is
\begin{displaymath}
  (p_{n+1},x_{n+1})=(a\cos(\omega t)+b\sin(\omega t),x_\mathrm{asy})
\end{displaymath}
where $x_\mathrm{asy}$ is the asymptotic output at the end of the
iteration, and $a$ and $b$ are the harmonic Fourier coefficients of the
asymptotic input $p(t)+u(t)$. The coefficients $a$ and $b$ are given
in \eqref{eq:pupdate} but they can also be extracted directly from the
input (the input $p(t)+u(t)$ is harmonic up to tolerance after the
iteration).

\section{Discussion of the method}
\label{sec:discussion}
The technique presented in section~\ref{sec:method:eq} is guaranteed
to be applicable near folds of equilibria involving one stable branch
for single inputs $p+u$ and outputs $x$ that give a bifurcation
diagram as in \fref{fig:simple-evolution}. It fails in points where
the equilibrium output $x$ does not depend on the parameter at the
linear level (that is, for example, at transcritical bifurcations or
when the branch is horizontal in the $(p,x)$-plane). The reason is
that the genericity assumption for successful control
(stabilizability and observability) is violated in these points.

As the control is applied by varying the bifurcation parameter $p$
several parameters are not truly independent. For example, in
section~\ref{sec:method:eq} the parameter pair $(\tilde p,\tilde x)$
enters only as a combination $\tilde p+k\tilde x$, such that one can
treat them as a single parameter. The same applies to the parameters
$a^*$ and $A^*_1$, and $b^*$ and $B^*_1$ in
section~\ref{sec:method:po}. Because a continuation of a forced system
in the forcing amplitude has effectively two free parameters (the
forcing amplitude and the phase), only the overall amplitude
$\sqrt{(a^*)^2+(b^*)^2+(A^*_1)^2+(B^*_1)^2}$ needs to be monitored.

The approach presented in sections \ref{sec:method:eq} and
\ref{sec:method:po} should be compared with the alternatives for
control-based identification of bifurcations suggested in the
literature: the detection of the fold bifurcation in \cite{ASFKRK99}
required identification of the normal form coefficients. Time-delayed
feedback and wash-out filtered feedback are not able to stabilize
equilibria uniformly near the fold \cite{HS05,H11,AWC94} (also when
they are modified by adding degrees of freedom \cite{TMPP07}). The
general approach proposed in \cite{SK08b} and taken in
\cite{SGNWK08,Barton2012,BB11,Bureau2011,Bureau2012} searches not for
the equilibrium on the branch that intersects the line given by the
feedback law \eqref{eq:controlinp}, but finds the equilibrium on a
prescribed line perpendicular to the secant through $(\tilde p,\tilde
x)$ (pseudo-arclength continuation). That is, it sets the feedback law
to
\begin{equation}
  \label{eq:controlinp2}
  u(t)=k(x^*-x(t))\mbox{,}
\end{equation}
and determines $(p^*,x^*)$ by solving the nonlinear system of
equations
\begin{align}
  \label{eq:fixpoint}
  x^*&=x_\mathrm{asy}(p^*,x^*)\mbox{,}\\
  \label{eq:pseudoarc}
  0&=\begin{bmatrix} p^*-\tilde p\\ x^*-\tilde x
  \end{bmatrix}^T
  \begin{bmatrix}
    p_n-p_{n-1}\\
    x_n-x_{n-1}
  \end{bmatrix}\mbox{,}
\end{align}
where $x_\mathrm{asy}$ is the steady-state output of the experiment
with parameter $p^*$ and feedback control \eqref{eq:controlinp2} after
the transients have settled.  This pseudo-arclength continuation in
\cite{SGNWK08,Barton2012,BB11,Bureau2011,Bureau2012} required
adjustments of $(p^*,x^*)$ in a Newton iteration for system
\eqref{eq:fixpoint},\,\eqref{eq:pseudoarc} to make the control
non-invasive (enforced by \eqref{eq:fixpoint}) on the line prescribed
by \eqref{eq:pseudoarc}. In this sense the procedure of
\fref{fig:simple-evolution} is a simplification of the control-based
pseudo-arclength continuation of
\cite{SGNWK08,Barton2012,BB11,Bureau2011,Bureau2012} that can be used
whenever the feedback control is applied by varying the bifurcation
parameter (which is not the case for \cite{Bureau2011,Bureau2012}).

The method presented in section~\ref{sec:method} promises a
substantial speed-up compared to
\cite{SGNWK08,Barton2012,BB11,Bureau2011,Bureau2012} by removing one
equation per free bifurcation parameter from the fixed point problem
\eqref{eq:fixpoint}, instead of adding the equation
\eqref{eq:pseudoarc}. The projection onto the solution surface occurs
along a line determined by the control gains $k$. Whenever the
remaining equations of \eqref{eq:fixpoint} can be solved by a simple
fixed-point iteration (or there are no equations remaining), this
removes the need for a full Newton iteration.

The extension to periodic orbits of forced systems proposed in
section~\ref{sec:method:po} applies the method shown in
\fref{fig:simple-evolution} to the harmonic part and combines it with
a simple fixed-point iteration for the non-harmonic part. In general,
the simple fixed-point iteration cannot be guaranteed to converge for
strongly non-harmonic periodic orbits. If one introduces a relaxation
parameter $R$ into the iteration \eqref{eq:picardj} (that is, setting
$X^*_{k+1}=(1-R)X^*_k+RX(X^*_k)$) and chooses $R$ small, then the
iteration becomes equivalent to the extended time-delayed feedback
(ETDF) method \cite{GSCS94} to finding periodic orbits (which also has
a relaxation parameter), but restricted to the non-harmonic Fourier
coefficients. This restriction to the non-harmonic Fourier
coefficients is essential. An algorithm updating all Fourier
coefficients in step 3 of the description in
section~\ref{sec:method:po} suffers from the same odd-number
limitation as the ETDF method.


\section{Experimental set-up and methods}
\label{sec:setup}

\begin{figure}[b]
  \centering
  \includegraphics{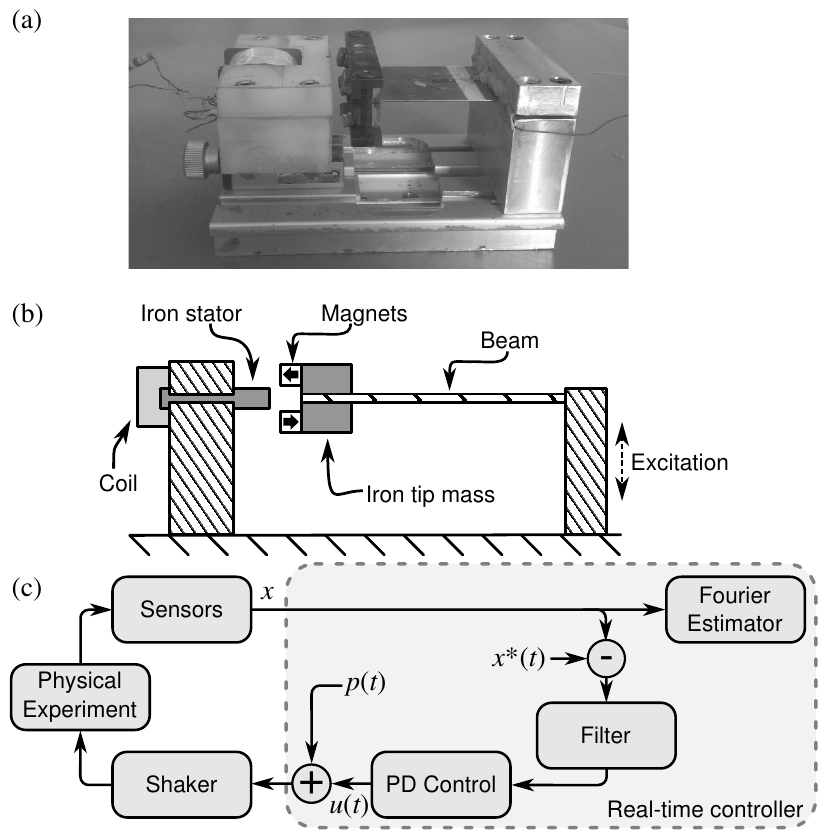}
  \caption{(a) A photograph of the nonlinear energy harvester. The
    approximate dimensions of the energy harvester are
    $137\,\text{mm}\times72\,\text{mm}\times85\,\text{mm}$ (width
    $\times$ height $\times$ depth). (b) A schematic of the physical
    parts of the nonlinear energy harvester. (c) A schematic of the
    experimental set-up. The elements within the shaded box are
    implemented within a real-time control system. The input $p(t)$
    are the system parameters, which in our case is harmonic forcing
    of the type $p(t)=a\cos(\omega t) + b\sin(\omega t)$.}
  \label{fig:setup}
\end{figure}
The demonstration experiment is the forced nonlinear oscillator shown
in \fref{fig:setup}, an electro-magnetic energy harvester
\cite{Barton2012,BB11}, mounted on a force-controlled electro-dynamic
shaker also shown in \fref{fig:setup}. 
Due to magnetic hysteresis and eddy currents this system is difficult
to characterize with a low number of degrees of freedom in a way that
is able to reproduce the experimental bifurcation diagram (see
\fref{fig:results}(a)) quantitatively. In particular, one would have
to introduce an ``effective'' damping coefficient that depends on the
forcing frequency and the response amplitude and phase
\cite{Cammarano2012}.

The experimental set-up, shown in \fref{fig:setup}, consists of a
generic electrodynamic shaker, a Maxon ADS 50/10 current controller
and a dSpace DS1104 real-time measurement and control
system. 

The input to the energy harvester is a force that is directly
proportional to the current supplied to the shaker. Thus the current
controller enables the force input to be determined directly. The
current-force relationship was determined by a series of quasi-static
tests.

The output $x$ is the displacement, measured from the energy harvester
using a suitably calibrated strain gauge.  The real-time controller
implements a fourth order IIR Butterworth filter ($-3$\,dB cut-off at
$75$\,Hz) and a proportional-derivative (PD) controller. The
derivative is estimated on-line using a two point finite difference
(the filter is sufficient to reduce the noise to an acceptable level
for a simple finite difference to work). The filter is purely for the
purposes of control; all other calculations use the unfiltered
data. 

The first seven Fourier coefficients ($m=7$) are estimated in
real-time from the unfiltered data using a recursive estimator to
minimize sampling and noise effects
caused by the forcing period not being an integer multiple of the
sampling period. While it is not necessary to calculate the
coefficients in real-time, it simplifies the implementation as it
reduces the communication needed between the real-time processor and
the host computer. The recursive estimator for the $k$-th Fourier
coefficient is
\begin{multline*}
  [A_{k,j+1},B_{k,j+1}] = \\
  [A_{k, j}, B_{k, j}] + \frac{\pi}{\omega}
  \int_{t-2\pi/\omega}^t[\cos(k\omega s),\sin(k\omega s)]\cdot\\
      [x(s) - A_{k,j}\cos(k\omega s) - B_{k,j}\sin(k\omega s)] \d s
\end{multline*}
where $A_{k,0} = B_{k,0} = 0$. A good approximation to the Fourier
coefficients is typically obtained within two iterations.

For the PD control $k_px+k_d\dot x$, the control gains $k_p=0.2$ and
$k_d=-0.004$ are kept constant throughout. For these control gains the
PD controller is globally stabilizing. This simplifies the methodology
presented in Sec.~\ref{sec:method:po} further: there is no need to
update the parameters $a$ and $b$ in the forcing term while the
experiment is running. That is, the input is of the form 
\begin{equation}\label{eq:pdinput}
  u(t)=k_p(x^*(t)-x(t)+k_d(\dot x^*-x(t))
\end{equation}
(no $p^*(t)$ here in contrast to
\eqref{eq:controlinppd}). The iteration of section~\ref{sec:method:po}
will reduce all non-harmonic coefficients of $u$ (after transients
have settled) below tolerance in a single step. With these
simplifications the protocol for tracking a branch of periodic orbits
in the amplitude is as follows. Denote the vectors of Fourier
coefficients as
\begin{align*}
  X^{\phantom{*}}&=[A_0,(A_j,B_j)_{j=1}^m] && \mbox{for the output $x(t)$,}\\
  X^*&=[A^*_0,(A^*_j,B^*_j)_{j=1}^m] && \mbox{for the reference $x^*(t)$,}\\
  U^{\phantom{*}}&=[A^u_0,(A^u_j,B^u_j)_{j=1}^m] && \mbox{for the input $u(t)$.}
\end{align*}
\begin{compactenum}
\item \label{step:picard0gen} Set $X^* := \tilde{X}_{n+1} = X_n +
  h[X_n-X_{n-1}]$ ($X_{n-1}$ and $X_n$ are the Fourier
  coefficients of outputs for the previous two points along the
  branch).
\item \label{step:controliterate1} Run the experiment with
  \eqref{eq:pdinput} and reference $x^*(t)$ corresponding to $X^*$
  until transients have settled. Then record the Fourier coefficients
  $X$ of the output $x(t)$.
\item \label{step:picard1} Set $X^* := X$ for all Fourier modes except
  the first ($A^*_1$ and $B^*_1$ are left unchanged).
\item \label{step:controliterate2} Run the experiment with
  \eqref{eq:pdinput} and reference $x^*(t)$ corresponding to $X^*$
  until transients have settled. Then record the Fourier coefficients
  $X$ and $U$ of the output and the control input respectively.
\end{compactenum}
The next point on the branch is then
\begin{align*}
  X_{n+1} &:= X\mbox{,}\\
  (a_{n+1},b_{n+1})&:=(A^u_1,B^u_1)
\end{align*}
(where $A^u_1$ and $B^u_1$ were recorded as part of $U$). All other
components of $U$ are zero up to experimental accuracy such that the
input $u(t)$ is indeed a harmonic forcing.

\paragraph*{Remarks}
\begin{compactitem}
\item The Fourier decomposition of $u$ and $x$ does not need to
  be done in real-time and instead can be done as a post-processing
  step to choose the new control target $X^*$ and to check
  convergence.
\item We accept the output as stationary when their corresponding
  Fourier coefficients become stationary for $5$ consecutive forcing
  periods.
\item The experiment is run continuously. That is, steps
  \ref{step:controliterate1} and \ref{step:controliterate2} of the
  procedure do not require initialization but continue from the state
  after previous steps.
\end{compactitem}
\section{Experimental results and discussion}
\label{sec:results}
\noindent 
We define three data measures, the forcing amplitude $F[u]$, the
response amplitude $R[x]$ and the error $e[u]$;
\begin{align}
  F[u]& := \sqrt{(A_1^u)^2 + (B_1^u)^2},\\
  R[x]& := \sqrt{(A_1)^2 + (B_1)^2},\\
  e[u]& := u(t) - A_1^u\cos(\omega t) -
  B_1^u\sin(\omega t),\label{eq:edef}
\end{align}
where the Fourier coefficients of the control input
$U=[A^u_0,(A^u_j,B^u_j)_{j=1}^m]$ and system response
$X=[A_0,(A_j,B_j)_{j=1}^m]$ are estimated continuously. When accepting
an output as a natural periodic orbit the error $e[u](t)$ should be
identically zero (to experimental accuracy).

\begin{figure}
  \centering
  \includegraphics{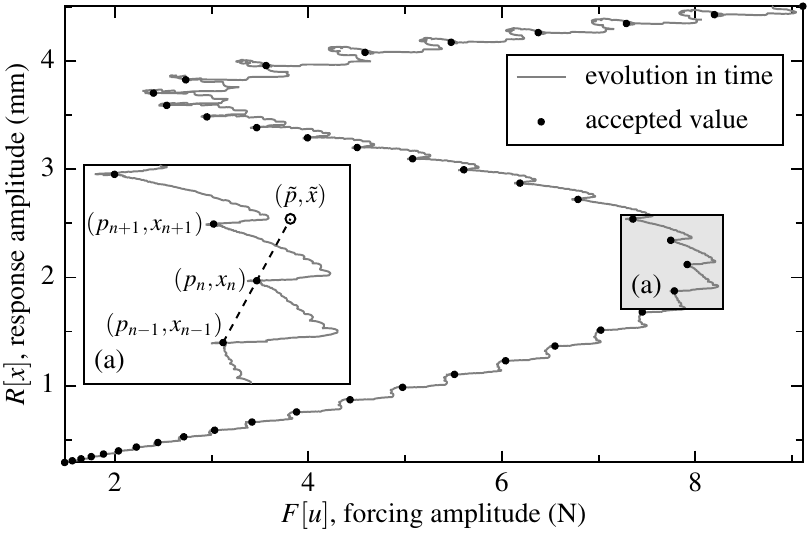}
  \caption{Experimental data showing the evolution of the controlled
    system as the bifurcation diagram of the uncontrolled system
    (shown in \fref{fig:setup}(a)). A family of periodic orbits is
    tracked through two saddle-node bifurcations (folds).
    Artificially large steps are taken along the solution curve for
    illustration purposes. Forcing frequency: $22$\,Hz.}
  \label{fig:evol:bif}
\end{figure}
\Fref{fig:evol:bif} shows the results of applying the methodology
described in Sec.~\ref{sec:setup} to the nonlinear energy
harvester. The tracking of periodic orbits starts from a stable,
low-amplitude periodic orbit and the forcing amplitude is then
increased.  As with \fref{fig:simple-evolution}, in
\fref{fig:evol:bif} two phases of the transients are visible between
the black dots: the horizontal coordinate $F[u]$ increases sharply
initially. This sharp increase is due to instantaneous changes in the
control target $x^*$. The rapid initial transient is followed by a
gradual stabilization towards the periodic orbit. Note that the output
measures $F[u]$ and $R[x]$ are not restricted to a single line but to
a higher-dimensional manifold because \fref{fig:evol:bif} is a
projection.  Hence, the point $(\tilde p,\tilde x)$ does not lie on
the line traced out by the evolution (in contrast to the sketch in
\fref{fig:simple-evolution}).

\begin{figure}
  \centering
  \includegraphics{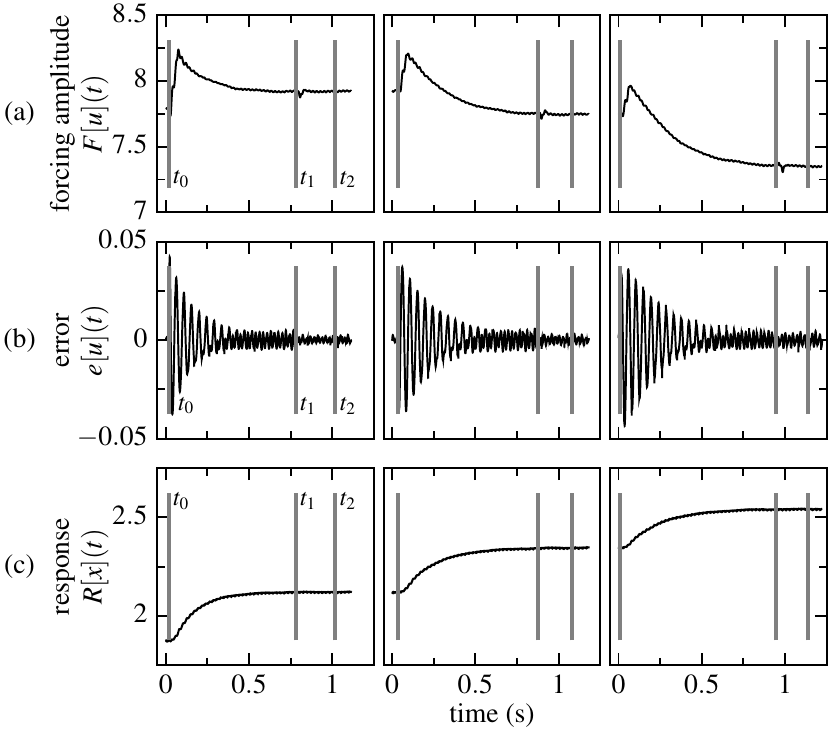}
  \caption{Time profile of forcing and response amplitudes and error
    (the non-harmonic part of control input $u$). Snapshots are time
    profiles corresponding to inset (a) of \fref{fig:evol:bif}. Note
    that the time gaps between the time profiles are only gaps in the
    time series recordings due to the saving of data (typically of the
    order of milli-seconds); the experiment ran continuously. Forcing
    frequency: $22$\,Hz, sampling frequency: $5$\,kHz.}
  \label{fig:timeseries}
\end{figure}
\begin{figure*}[t]
  \parbox{11cm}{
    \includegraphics{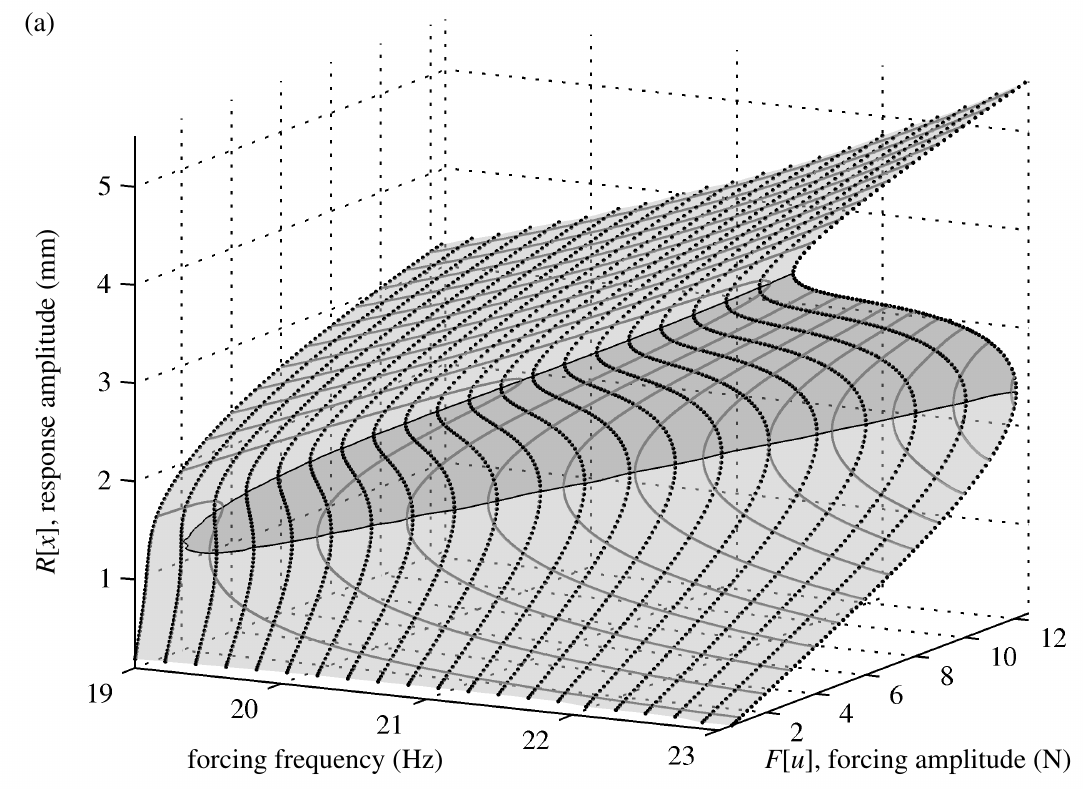}
  }
  \hfill
  \parbox{6.3cm}{
    \includegraphics{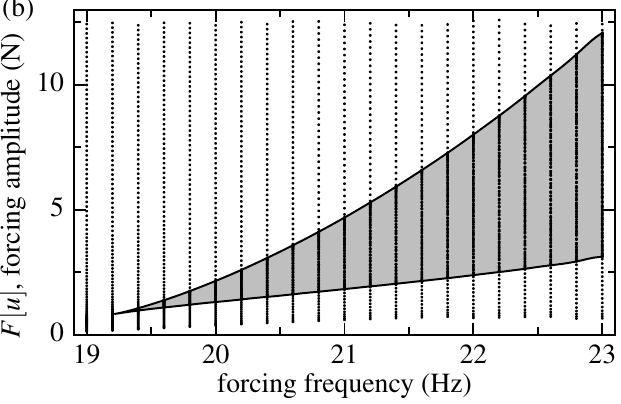}

    \includegraphics{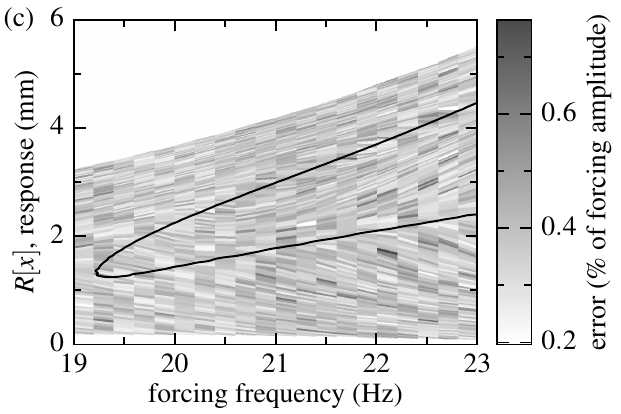}
  }
  \caption{Experimental results from the energy harvester shown in
    \fref{fig:setup}. A sequence of constant forcing frequency runs
    were performed at a spacing of $0.2$\,Hz. Panel (a) shows the
    complete resonance surface of the oscillator. Panel (b) shows the
    corresponding two-parameter bifurcation diagram (a top-down view
    of panel (a)) with the cusp point evident at approximately
    $19.2$\,Hz. Panel (c) is a front view of the resonance surface
    with the measured error superimposed onto the surface. The error
    is defined as the root-mean-square (RMS) of $e[u]$ as a percentage
    of the forcing amplitude; it measures how invasive the method
    is. In all panels, the data points are shown as black dots and the
    calculated saddle-node bifurcation (fold) curve is marked in
    black. Points within the dark gray region of panels (a) and (b)
    are unstable solutions.}
  \label{fig:results}
\end{figure*}
\Fref{fig:timeseries} shows the time series recordings corresponding
to the data points shown in \fref{fig:evol:bif}(a); they demonstrate
in detail the convergence of the method as the system passes through a
saddle-node bifurcation (fold).  The first harmonic of the input
$F[u]$ (\fref{fig:timeseries}(a)) gradually drifts during non-periodic
transients but settles rapidly. The error $e[u]$
(\fref{fig:timeseries}(b)) corresponds to the non-harmonic, invasive,
part of the control; its decay during each step is evident. The
vertical bars (marked $t_0$, $t_1$ and $t_2$) indicate the stages of
the iteration: step~\ref{step:controliterate1} occurs from $t_0$ to
$t_1$, and step~\ref{step:controliterate2} occurs from $t_1$ to
$t_2$. The input $u$ and output $x$ after $t_2$ are then accepted as
points on the branch. \Fref{fig:timeseries}(c) shows the amplitude of
the first harmonic $R[x]$ of the displacement to further demonstrate
convergence.

The main advantage of the method presented here over methods based on
Newton iterations, apart from ease of implementation, is the speed-up
of a factor of $\approx15$ compared to \cite{SK08b,Barton2012} (a
conservative estimate; only individual solution curves could be traced
out in \cite{SK08b,Barton2012}). This feature is particularly
important if one wants to explore systems that gradually degrade under
laboratory conditions.

As illustrated in \fref{fig:results}(a), this speed-up enables
tracking of entire surfaces and the associated bifurcations. The
experimental data points (marked by black dots in panels~(a) and (b))
are obtained by consecutive runs for fixed frequencies $0.2$\,Hz
apart. The total experimental time to generate these results was 61
minutes. Panel~(a) shows the three-dimensional projection in the space
spanned by the two parameters forcing frequency and amplitude and the
response amplitude (note that the response is non-harmonic and so this
is indeed only a projection). Its main feature is the curve of
saddle-node bifurcations passing through a cusp bifurcation
(black). To facilitate the extraction of geometric information, the
data points in \fref{fig:results} are interpolated using Wendland's
compactly supported radial basis
functions~\cite[Ch.~11]{Fasshauer2007}. Using the interpolant, the
bifurcation and constant forcing amplitude curves in
\fref{fig:results}(a,b) are calculated using numerical continuation on
the experimentally generated surface.  Curves of constant forcing
amplitude (gray), reminiscent of the resonance curves for an idealized
Duffing oscillator, give additional geometric information.
All the data points within the dark gray shaded area of
\fref{fig:results}(a) are unstable periodic orbits of the uncontrolled
system, and would typically not be seen experimentally.

\Fref{fig:results}(b) shows a top-down view of panel (a), a
two-parameter bifurcation diagram, again indicating all measured
points on the unstable part of the surface in a darker shade
of gray. The saddle-node bifurcation (fold) curve bounds the
instability region with a cusp point at approximately
$19.2$\,Hz.

\Fref{fig:results}(c) shows a front view of panel (a) with the error
at each data point rendered onto the surface. Here the error is
defined as the root-mean-square (RMS) of the non-harmonic part $e[u]$
(defined in \eqref{eq:edef}) over one period. This is a
measure of the invasiveness of the control; if this method was truly
non-invasive, then this error would be zero. In the experimental
set-up here the error is low, with a mean error of $<0.5\%$. The
predominant source of error is noise amplification through the use of
a derivative controller. This error is kept to a minimum through the
use of the Butterworth filter described in section~\ref{sec:setup}. As
seen in \fref{fig:results}(c), there is no apparent correlation
between geometric features of the solution surface (e.g., the fold
points) and the magnitude of the error at that point.

\section{Conclusion}
\label{sec:conclusion}

The presented approach is a general experimental technique to explore
dynamical systems in parameter studies near saddle-node bifurcations
(folds). It is particularly useful for the exploration of families of
equilibria because no iterations similar to \eqref{eq:picardj} are
necessary. As we demonstrated, it is also applicable to periodically
forced systems (the generalization to a non-harmonic forcing is
straightforward). The main limitation of the method is that control
fails at the linear level whenever the system does not depend on the
bifurcation parameter to first order (e.g., near transcritical
bifurcations).  As the presented approach works particularly well
around saddle-nodes, its main application areas are likely
complementary to those of Pyragas' TDF control. Examples currently
under investigation include the identification of growth rates in
chemostats \cite{VFDNLS12,RSRD12}, or tracking localized spots in
ferrofluids \cite{GRR10}.

There are several ways in which this method can be generalized. First,
if the equilibrium has more than a single unstable dimension, one
typically reconstructs a proxy for the state through an observer
$x_\mathrm{obs}$ \cite{S98} and lets $u$ depend on
$x_\mathrm{obs}$. This is a generalization of the use of PD control in
the sections \ref{sec:method:po} and \ref{sec:setup}. Second, if one
has more than a single adjustable system parameter then one can obtain
multi-parameter families of equilibria and apply feedback control
through more than a single input. Similarly, if one uses more than a
single output ($x\in\R^m$, $m>1$), one can feed back the input $u$
depending on the multi-dimensional $x$, making control easier to
achieve.

Acknowledgments: The research of J.S.\ is supported by EPSRC Grant
EP/J010820/1.

\bibliography{expbif}
\end{document}